\documentclass{amsart}

\usepackage{amssymb,amsfonts, latexsym,amsmath,hhline,array,longtable}
\usepackage{pdfsync,color, comment,colortbl, mathrsfs,stmaryrd,cite,graphicx,yfonts}

\usepackage{ifpdf}
\ifpdf \usepackage[colorlinks=true, citecolor=blue, linkcolor=blue, urlcolor=blue]{hyperref} \fi

\newcommand{\cal}{\mathcal}

\newtheorem{formula}{}[section]
\newtheorem{definition}[formula]{Definition}
\newtheorem{corollary}[formula]{Corollary}
\newtheorem{remark}[formula]{Remark}
\newtheorem{lemma}[formula]{Lemma}
\newtheorem{theorem}[formula]{Theorem}
\newtheorem*{claim}{Claim}

\def\thrm{\begin{theorem}}
\def\thrml#1{\begin{theorem}\label{#1}}
\def\ethrm{\end{theorem}}
\def\rmrk{\begin{remark}}
\def\rmrkl#1{\begin{remark}\label{#1}}
\def\ermrk{\end{remark}}
\def\dfntn{\begin{definition}}
\def\dfntnl#1{\begin{definition}\label{#1}}
\def\edfntn{\end{definition}}
\def\nmrt{\begin{enumerate}}
\def\enmrt{\end{enumerate}}
\def\tm#1{\item[{\rm (#1)}]}
\def\qtnl#1{\begin{equation}\label{#1}}
\def\eqtn{\end{equation}}
\def\lmm{\begin{lemma}}
\def\lmml#1{\begin{lemma}\label{#1}}
\def\elmm{\end{lemma}}
\def\crllr{\begin{corollary}}
\def\crllrl#1{\begin{corollary}\label{#1}}
\def\ecrllr{\end{corollary}}
\def\css{\begin{cases}}
\def\ecss{\end{cases}}
\def\prf{\begin{proof}}
\def\eprf{\end{proof}}
\def\clm{\begin{claim}}
\def\eclm{\end{claim}}

\def\cA{{\cal A}}

\def\cP{{\cal P}}

\def\cS{{\cal S}}

\def\mC{{\mathbb C}}
\def\mF{{\mathbb F}}
\def\mN{{\mathbb N}}

\def\mZ{{\mathbb Z}}

\def\fK{{\mathfrak K}}
\def\fL{{\mathfrak L}}

\def\fX{{\mathfrak X}}

\DeclareMathOperator{\aut}{Aut}

\DeclareMathOperator{\Span}{Span}

\DeclareMathOperator{\sym}{Sym}

\def\qaq{\quad\text{and}\quad}

\def\ov{\overline}
\def\und#1{{\underline{#1}}}

\newtheorem{proposition}[formula]{Proposition}
\newtheorem{example}[formula]{Example}
\newtheorem*{introtheorem}{Theorem}
\numberwithin{equation}{section}

\begin{document}

\title{Multivalued Groups and Schur Rings}
\author{A.\,V.~Vasil'ev}
\address{Sobolev Institute of Mathematics, Siberian Branch of the Russian
Academy of Sciences, 4 Akademika Koptyuga Avenue, Novosibirsk 630090, Russia}
\email{vasand@math.nsc.ru}
\author{I.\,N.~Ponomarenko}
\address{St. Petersburg Department of Steklov Mathematical Institute of the
Russian Academy of Sciences, 27 Fontanka Embankment, St. Petersburg 191023, Russia}
\email{inp@pdmi.ras.ru}
\date{}

\begin{abstract}
We study finite involutive multivalued groups associated with Schur
rings. Their reduced valency is shown to equal the least common multiple
of the sizes of the basic sets of the corresponding ring. We prove that
the double-coset group associated with a complemented subgroup of a
finite group is isomorphic to a Schur multivalued group over the
complement. In particular, every coset multivalued group is a Schur
multivalued group. Using affine Schur rings, we construct families of
pairwise nonisomorphic commutative multivalued groups over elementary
abelian groups. We obtain asymptotic bounds for the sizes of these
families and prove that the proportion of double-coset groups in them
tends to zero.
\smallskip

\end{abstract}

\maketitle

\section{Introduction}
The theory of multivalued groups goes back to a construction of Buchstaber
and Novikov~\cite{BuchN1971}; a survey of its main results and applications
is given in~\cite{Bu2006}. Natural examples of multivalued groups arise in
topology~\cite{BukH1982}, algebra and dynamics~\cite{BuV2019,BukK2026}, and
algebraic combinatorics~\cite{Yagod2002,GPV2025}. Two important sources of examples
are the coset and double-coset constructions. In the former, a finite group
acts by automorphisms on another finite group, and the elements of the
resulting multivalued group are the orbits of this action. In the latter,
the elements of the multivalued group are the double cosets of a subgroup
of a finite group. We denote the classes of coset and double-coset
multivalued groups by $\fK_{coset}$ and $\fK_{bicoset}$, respectively.

In this note, we study multivalued groups arising from Schur rings. This
construction is part of the general theory of combinatorial
algebras~\cite{BuVEP1996e}; here we clarify its relationship with double-coset
groups and construct a family in which the proportion of double-coset
groups tends to zero. A Schur ring over a finite group~$G$ can be defined
as a subspace $\cA$ of the group algebra $\mC G$, that contains the identities $1_G$ and
$\sum_{g\in G}g$ for ordinary and componentwise multiplications, respectively,
and is closed under both multiplications and the linear extension of the
map $g\mapsto g^{-1}$. A definition in terms of a partition of the group
is given in Section~3.

Wielandt developed the theory of Schur rings in connection with permutation
groups containing a regular subgroup; see~\cite[Chapter IV]{Wielandt1964}.
Schur rings are now used in the theory of Cayley graphs~\cite{KlinP1978},
algebraic combinatorics~\cite{MuzychukR2024}, algebraic coding
theory~\cite{WuYY2022}, algorithmic problems in group and graph
theory~\cite{PonomV2017}, and other areas of mathematics (see, in particular, the survey \cite{MP2009}).

In Theorem~\ref{151125a}, we apply the construction
of~\cite{BuVEP1996e} to Schur rings. The normalized sums of their basic sets
form a basis of the group algebra of an involutive multivalued group,
whose reduced valency is the least common multiple of the sizes of these
sets (Theorem~\ref{151125a}). A multivalued group isomorphic to one obtained
in this way from an S-ring over $G$ will be called a Schur multivalued
group over $G$.

Every coset multivalued group is a Schur multivalued group
(Corollary~\ref{110426t}). More generally, if a subgroup $L$ of a finite
group $K$ has a complement~$G$, then the double-coset group of the pair
$(K,L)$ is isomorphic to a Schur multivalued group over~$G$
(Theorem~\ref{020426a}). For the class $\fK_{\mathrm{schur}}$ of all Schur
multivalued groups, we have the inclusions
\qtnl{121026a}
\fK_{coset}\,\subsetneq\,\fK_{bicoset}\,\cap\,\fK_{schur}
\,\subsetneq\,\fK_{schur}.
\eqtn
The first inclusion is strict by Example~\ref{270326a}; the strictness of
the second follows from the theorem below. It provides families of Schur
multivalued groups of bounded order whose sizes grow subexponentially
and in which the proportion of double-coset groups tends to zero.

\begin{introtheorem}
For every prime power $q\ge3$, there exists a family $\fK(q)$ of pairwise
nonisomorphic involutive multivalued groups with the following properties:
\nmrt
\tm{1} Each group in $\fK(q)$ has order at most $\lfloor(q+3)/2\rfloor$
and reduced valency at least $q-1$, and is a Schur multivalued group over
an elementary abelian group of order $q^2$.
\tm{2} $|\fK(q)|=\exp(\Theta(\sqrt q))$ as $q\to\infty$.
\tm{3} The proportion of double-coset groups in $\fK(q)$ tends to $0$
as $q\to\infty$.
\enmrt
\end{introtheorem}

Section~2 collects the basic definitions and properties of multivalued
groups. Schur multivalued groups are introduced in Section~3, and their
relationship with double-coset groups is established in Section~4. The theorem
in the introduction follows from Theorems~\ref{120426h} and~\ref{180426a}
and the estimate for the number of partitions given in the proof of the
latter.

\section{Multivalued groups}\label{170824d}

\subsection{Definition}
Let $X$ be a finite nonempty set, let $n$ be a positive integer, and let
$X^{(n)}$ be the set of all $n$-multisets of elements of $X$. An $n$-valued
group structure on $X$ is specified by a multiplication
$\cdot:X\times X\to X^{(n)}$ and an inversion map $\star:X\to X$
satisfying the conditions below. We write the product of $x,y\in X$ as
\[
x\cdot y=[(x\cdot y)_1,\ldots,(x\cdot y)_n].
\]
Multiplication extends to multisets in the natural way. It is required
to be associative: for all $x,y,z\in X$, the following equality of
$n^2$-multisets holds:
\[
(x\cdot y)\cdot z=x\cdot(y\cdot z).
\]
We also assume that there is an identity element $e\in X$ for
multiplication, that is,
\[
x\cdot e=e\cdot x=[\underbrace{x,\dots,x}_{n}]
\]
for all $x\in X$. The inverse $x^\star$ satisfies $e\in x\cdot x^\star$
and $e\in x^\star\cdot x$. The number $|X|$ is called the order of the
multivalued group~$\fX$.

For $x,y,z\in X$, let $m_{x,y}^z$ denote the multiplicity of $z$ in the
multiset $x\cdot y$. Setting $p_{x,y}^z=m_{x,y}^z/n$, we obtain the group
algebra of the multivalued group $\fX$: it is the complex vector space
with basis $X$ and multiplication $xy=\sum_zp_{x,y}^z z$. We call the number
\[
\nu(\fX)=\min\{N\in\mN: Np_{x,y}^z\in\mZ\text{ for all }x,y,z\in X\}
\]
the \emph{reduced valency}. Thus, $\nu(\fX)$ is the least common multiple
of the positive denominators of the rational numbers
$p_{x,y}^z$ written in lowest terms.

Following~\cite{BuVEP1996e}, we use a normalized notion of isomorphism of multivalued groups
that allows the valencies to differ. Namely, a bijection
$f:X\to X'$, $x\mapsto x'$, is an \emph{isomorphism} from an $n$-valued
group $\fX$ to an $n'$-valued group $\fX'$ with underlying set~$X'$ if
\qtnl{110426f}
(x^\star)'=(x')^\star\qaq p_{x^{},y^{}}^{z^{}}=p_{x',y'}^{z'}
\eqtn
for all $x,y,z\in X$ (in particular, $f(e)=e'$). If $n=n'$, the latter
equality is equivalent to the usual preservation of multiset
multiplication. The valency $n$ is not an invariant of isomorphism as
defined here, whereas the reduced valency is.

In accordance with~\cite[Definition 1.2]{BuVEP1996e}, a multivalued group
$\fX$ is called \emph{involutive} if
\[
(x^\star)^\star=x,\quad e\in x\cdot y\Longleftrightarrow y=x^\star,\quad
m_{x,x^\star}^e=m_{x^\star,x}^e,\quad (x\cdot y)^\star=y^\star\cdot x^\star.
\]
The first and last conditions mean that $\star$ is an involutive
anti-isomorphism of the group algebra of~$\fX$.

A subset $Y\subseteq X$ is called a \emph{subgroup} of the multivalued
group $\fX$ if $e\in Y$, $Y^\star=Y$, and $x\cdot y$ consists of elements
of $Y$ for all $x,y\in Y$. The restricted operations give $Y$ an
$n$-valued group structure. 

\subsection{Double-coset groups}\label{160426a}
Let $K$ be a finite group and $L$ a subgroup of order $n$.
Following~\cite[Section~6]{Bu2006}, let $X=L\backslash K\slash L$ be the set of
double cosets, and define an $n$-valued multiplication by
\[
m_{x,y}^z=|\{\ell\in L: Lk_1\ell k_2L=z\}|,
\]
where $x=Lk_1L$ and $y=Lk_2L$. With this multiplication, $X$ becomes an
$n$-valued group $\fX$ with identity element $L$ and inversion
$(LkL)^\star=Lk^{-1}L$, $k\in K$. We call $\fX$ the \emph{double-coset}
group associated with the pair $(K,L)$. Note that $\star$ is an involution
on $X$ and $(x\cdot y)^\star=y^\star\cdot x^\star$. The uniqueness of the
inverse and the equalities
\[
m_{x,x^\star}^{L}=|L\cap k^{-1}Lk|,
\qquad m_{x^\star,x}^{L}=|L\cap kLk^{-1}|
\]
for $x=LkL$ imply that $\fX$ is involutive.

For the double-coset group of $(K,L)$, let $\ov K$ denote the permutation
group on the set $L\backslash K$ of right cosets induced by right multiplication in $K$. Note
that $\ov K$ is the image of the action of $K$ in the symmetric group;
this action need not be faithful.

\lmml{160426h}
The permutation group $\ov K$ has degree
$\sum_{x\in X}(p_{x^\star,x}^e)^{-1}$. It is primitive\footnote{Here we
also regard the group of degree~$1$ as primitive.} if and only if the
multivalued group $\fX$ has no nontrivial proper subgroups.
\elmm
\prf
The degree of $\ov K$ is $|K:L|=\sum_{x\in X}c(x)$, where $c(x)$ is the
number of right cosets of $L$ contained in $x$. Let $x=LkL$. Then
$c(x)=|LkL|/|L|$. The standard formula for the size of a double coset gives
\[
n\,c(x)=|LkL|=\frac{|L|\cdot|L|}{|L\cap kLk^{-1}|}
=\frac{n^2}{m_{x^\star,x}^e}=\frac{n}{p_{x^\star,x}^e}.
\]
Thus $c(x)=(p_{x^\star,x}^e)^{-1}$, as required.

Let us prove the second assertion of the lemma. It is clear if $L=K$. Suppose that $L<K$. The action of $K$
on $L\backslash K$ is primitive exactly when $L$ is a maximal subgroup
of $K$. The subgroups $H$ with $L\le H\le K$ are in bijection with the
subgroups of the double-coset group: to $H$ we associate the set of double
cosets contained in $H$, and to a subgroup $Y$ of $\fX$ we associate the union
$\bigcup_{x\in Y}x$. Hence the absence of intermediate subgroups between
$L$ and $K$ is equivalent to the absence of subgroups $Y$ with
$\{L\}\subsetneq Y\subsetneq X$.
\eprf

We conclude this section with an important special case of double-coset
groups. Let $G$ and $A$ be finite groups and let
$\varphi:A\to\aut(G)$ be a homomorphism. Set $K=G\rtimes_\varphi A$
and $L=\{1\}\rtimes A$. Identifying $L\backslash K$ with $G$, we see that
the right action of $K$ on $L\backslash K$ induces on $G$ the action
generated by right translations and the automorphisms in $\varphi(A)$.
The double cosets in $L\backslash K\slash L$ correspond to the
$\varphi(A)$-orbits on $G$. Following \cite{Bu2006}, we call the resulting multivalued group on the
set of $\varphi(A)$-orbits the \emph{coset} group of the triple
$(G,A,\varphi)$.

\section{Schur multivalued groups}
\subsection{Schur rings}
We recall the definition of an S-ring. Let $G$ be a finite group. A
subalgebra $\cA$ of the group algebra $\mC G$ is called a \emph{Schur ring},
or an \emph{S-ring}, over~$G$ (see~\cite{Sch1933,Wielandt1964}) if there
exists a partition $\cS=\cS(\cA)$ of $G$ such that $\{1_G\}\in\cS$,
$R^{-1}\in\cS$ for every $R\in\cS$, and
\qtnl{210326a}
\cA=\Span_{\mC}\{\und{R}: R\in\cS\},
\eqtn
where $\und{R}=\sum_{g\in R}g$. The classes of the partition $\cS$ are
called the basic sets of~$\cA$.

Note that $\cS$ is uniquely determined by the algebra~$\cA$. The largest
S-ring over $G$ corresponds to the discrete partition of $G$ and is equal
to $\mC G$. If $G\ne1$, the smallest S-ring corresponds to the partition
with two classes $\{1_G\}$ and $G^\#=G\setminus\{1_G\}$; this S-ring is
called trivial.

\begin{example}\label{260326a}
{\rm\cite{Sch1933}; see also~\cite[Section~21]{Wielandt1964}.}
Let $K$ be a subgroup of the symmetric group $\sym(G)$ containing the
permutations $x\mapsto xg$, $x\in G$, for every $g\in G$. Let $\cS$ be
the partition of $G$ into the orbits of the stabilizer $K_1=K_{1_G}$ of
$1_G$. Then the subspace $\cA$ defined by~\eqref{210326a} is an S-ring
over $G$. We denote it by $\fL(G,K_1)$.
\end{example}

It follows from~\eqref{210326a} that the elements $\und{R}$, $R\in\cS$,
form a linear basis of $\cA$. In particular, the product of two such
elements is a linear combination
\qtnl{080426t}
\und{R}\,\und{S}=\sum_{T\in\cS}
c_{\scriptscriptstyle R,\scriptscriptstyle S}^{\scriptscriptstyle T}\und{T}
\eqtn
with nonnegative integer coefficients
$c_{\scriptscriptstyle R,\scriptscriptstyle S}^{\scriptscriptstyle T}$.
These coefficients are called the structure constants of the S-ring and
satisfy the identities
\qtnl{1210266b}
|T|c_{R,S}^{T}=|R|c_{T,S^{-1}}^{R}\qaq
\sum_T|T|c_{R,S}^{T}=|R||S|
\eqtn
for all $R,S,T\in\cS$. It is easy to check that
$c_{\scriptscriptstyle R,\scriptscriptstyle S}^{\scriptscriptstyle T}$
is the number of representations of a fixed $t\in T$ as $t=gh$, with
$g\in R$ and $h\in S$. Note that
\[
d(R)=c_{\scriptscriptstyle R,\scriptscriptstyle R^{-1}}^{1_G}=|R|
\]
is positive, and the sum of these numbers over all $R\in\cS$ equals the
order of $G$ (to simplify notation, here and below we write $1_G$ instead
of $\{1_G\}$ in the superscript). Moreover, $d(R)=d(R^{-1})$.

For each $A=\sum_{T\in\cS}a_T\und{T}$ in $\cA$, set
\[
A^*=\sum_{T\in\cS}\bar a_T\und{T^{-1}},
\]
where $\bar a_T$ is the complex conjugate of $a_T$. The bijection
$A\mapsto A^*$, $A\in\cA$, is a semilinear involutive antiautomorphism of
$\cA$: $(AB)^*=B^*A^*$ and $(A^*)^*=A$. The map
\[
\varepsilon:\cA\to\mC,\quad
\sum_Ta_T\und{T}\mapsto\sum_Ta_T|T|,
\]
is a one-dimensional representation and satisfies
$\varepsilon(A^*)=\ov{\varepsilon(A)}$.

\subsection{Construction}
We construct a multivalued group from an S-ring $\cA$. Set
\qtnl{250326a}
n=\operatorname{LCM}\{|R|:R\in\cS\}\qaq
X=\{x_R:R\in\cS\},
\eqtn
where $x_R=\frac{1}{|R|}\und{R}$. Define $n$-valued multiplication on
$X$ by prescribing the multiplicity $m_{x,y}^z$ of $z=x_T$ in the product
$x\cdot y$ of $x=x_R$ and $y=x_S$ as follows:
\qtnl{250326a1}
m_{x,y}^z=n\cdot
\frac{|T|\,c_{\scriptscriptstyle R,\scriptscriptstyle S}^{\scriptscriptstyle T}}
{|R|\,|S|}.
\eqtn
Using the identities~\eqref{1210266b}, one readily checks that $m_{x,y}^z$
is a nonnegative integer and $\sum_zm_{x,y}^z=n$. Finally, set
$(x_R)^\star=x_{R^{-1}}$.

\thrml{151125a}
With the notation above, the triple $\fX=\fX(\cA)=(X,\cdot,\star)$ is an involutive $n$-valued group,
and $\nu(\fX)=n$.
\ethrm
\prf
Set $p_{R,S}^T=m_{x,y}^z/n$. The established properties of $m_{x,y}^z$ and
the definition of~$n$ imply that $np_{R,S}^T$ is a nonnegative integer and
$\sum_Tp_{R,S}^T=1$. The associativity of multivalued multiplication
follows from the associativity of multiplication in $\cA$ and the
equalities
\[
x_Rx_S=\sum_{T\in\cS}p_{R,S}^T x_T.
\]
The identity is $x_{\{1_G\}}=1_G$. Since
\[
p_{R,S}^{1_G}=\css
1/|R|,&\text{if $S=R^{-1}$},\\
0,&\text{if $S\ne R^{-1}$},
\ecss
\]
the inverse is unique and equals $x_{R^{-1}}$. The equality
$|R|=|R^{-1}|$ gives equality of the multiplicities of the identity in
$x_R\cdot x_{R^{-1}}$ and $x_{R^{-1}}\cdot x_R$. The property
$(AB)^*=B^*A^*$ ensures that the involution is compatible with
multiplication. Thus the group is involutive.

All the numbers $np_{R,S}^T$ are integers, so $\nu(\fX)$ divides $n$.
On the other hand, the normalized multiplicities include
$p_{R,R^{-1}}^{1_G}=1/|R|$ for every $R\in\cS$. Hence each $|R|$ divides
$\nu(\fX)$, and therefore $n$ divides $\nu(\fX)$. Thus $\nu(\fX)=n$.
\eprf

\dfntn
We call $\fX(\cA)$ the multivalued group associated with the S-ring
$\cA$. A group $\fX$ is called a Schur multivalued group over $G$ if it
is isomorphic to $\fX(\cA)$ for some S-ring $\cA$ over $G$.
\edfntn

The following example shows that S-rings over different groups may give
rise to isomorphic Schur multivalued groups.

\begin{example}\label{270326a}
Let $\cA$ be the trivial S-ring over a group $G$ of order $m>1$.
Formula~\eqref{250326a} and Theorem~{\rm \ref{151125a}} show that
$\fX:=\fX(\cA)$ is an $(m-1)$-valued group of order~$2$. Its
multiplicities~\eqref{250326a1} depend only on $m$, not on $G$. Thus
$\fX$ is a Schur multivalued group over every group $G$ of order $m$.

Moreover, $\fX$ is a double-coset group: it is isomorphic to the group
associated with the pair $(\sym(m),\sym(m-1))$, where $\sym(m-1)$ is a
point stabilizer in the natural action of $\sym(m)$. Indeed, both
multivalued groups have two elements, and the normalized multiplicity
of the identity in the square of the nonidentity element is $1/(m-1)$.

If $m$ is divisible by two distinct primes $p$ and $q$, then $\fX$ is not
a coset group. Otherwise, it would arise from a transitive action of a
group of automorphisms on the nonidentity elements of a group $G$ of
order $m$. However, $G$ contains elements of orders $p$ and $q$, which
cannot belong to the same orbit of a group of automorphisms.
\end{example}

\subsection{Isomorphisms}
Isomorphic multivalued groups can also arise from different S-rings over
the same group. To describe this phenomenon, we recall the notion of an
algebraic isomorphism.

Let $\cA$ and $\cA'$ be S-rings over groups $G$ and $G'$, respectively.
An isomorphism of complex algebras $\varphi:\cA\to\cA'$ is called an
\emph{algebraic isomorphism}\footnote{Combinatorial isomorphisms of
S-rings are defined as isomorphisms of the corresponding Cayley schemes.}
if it maps the standard basis of $\cA$ onto the standard basis of $\cA'$.
In this case, for every $R\in\cS(\cA)$ there is a unique
$R'\in\cS(\cA')$ such that $\varphi(\und{R})=\und{R'}$. Moreover,
\[
\{1_G\}'=\{1_{G'}\},\quad (R^{-1})'=(R')^{-1},\quad
c_{\scriptscriptstyle R,\scriptscriptstyle S}^{\scriptscriptstyle T}
=c_{\scriptscriptstyle R',\scriptscriptstyle S'}^{\scriptscriptstyle T'}
\]
for all $R,S,T\in\cS(\cA)$. Therefore, \eqref{250326a}
and~\eqref{250326a1} show that an algebraic isomorphism~$\varphi$
induces an isomorphism $x_{R^{}}\mapsto x_{R'}$, $R\in\cS(\cA)$,
of the associated multivalued groups, thus proving the following statement.

\begin{proposition}\label{300326a}
If $\cA$ and $\cA'$ are algebraically isomorphic S-rings, then the
multivalued groups $\fX(\cA)$ and $\fX(\cA')$ are isomorphic and have
the same reduced valency.
\end{proposition}

In Section~\ref{300326b}, we construct a family of S-rings over an elementary
abelian group of order $q^2$ whose number of pairwise algebraically
nonisomorphic representatives grows as $\exp(\Theta(\sqrt q))$
when $q\to\infty$.

\section{Double-coset multivalued groups in the context of Schur rings}
Let $K$ be a finite group. The double cosets of a subgroup $L$ of $K$
form a partition of $K$ closed under taking inverses. The linear span of
their sums is closed under multiplication in $\mC K$. However, the double
coset containing the identity of $K$ is $L$, rather than $\{1_K\}$.
Thus, if $L\ne1_K$, this partition does not define an S-ring over $K$ in
the usual sense. According to the definition in~\cite[Section~3]{Muzychuk1991},
such a partition defines only a pseudo-S-ring over $K$. We show that,
under an additional assumption on $L$, the corresponding multivalued
group $\fX$ is associated with a suitable Schur ring.

Suppose that $L$ has a complement $G$ in $K$, that is,
\[
K=L\,G\qaq G\cap L=1.
\]
Then each $k\in K$ has a unique factorization $k=\ell g$, where $\ell\in L$
and $g\in G$. This defines a bijection
\[
\pi:L\backslash K\longrightarrow G,\qquad \pi(Lg)=g.
\]
Via this bijection, the permutation group $\ov K$ defined
in Section~\ref{160426a} can be identified with the group of all permutations
$\ov k\in\sym(G)$ determined by $k\in K$ through
\[
g^{\ov k}=\pi(Lgk),\quad g\in G,\ k\in K.
\]
Let $N$ be the kernel of the corresponding action of $K$ on $G$. Then
$\ov K_{1_G}\cong L/N$. Moreover, $\ov K\cong K/N$ contains the right
regular subgroup $G$. Thus Example~\ref{260326a} defines the S-ring
\[
\cA(K,L)=\fL(G,\ov K_{1_G}).
\]
In general, this S-ring depends on the choice of the complement $G$;
this choice is suppressed in the notation $\cA(K,L)$.

\thrml{020426a}
Let $L$ be a subgroup of a finite group $K$, and let $G$ be a complement
to $L$ in $K$. Then the double-coset group of the pair $(K,L)$ is
isomorphic to the multivalued group associated with the S-ring
$\cA(K,L)$ over $G$.
\ethrm
\prf
Set $X=L\backslash K\slash L$. For each $x\in X$, let
\[
\pi(x)=\{\pi(Lk):k\in x\},
\]
where $\pi$ is defined above. If $x=LgL$, then the right cosets of $L$ contained in $x$ are $Lg\ell$,
$\ell\in L$. Thus the sets $\pi(x)$, $x\in X$, are precisely the
$\ov K_{1_G}$-orbits, that is, the basic sets of $\cA(K,L)$. We verify
that for all $x,y,z\in X$,
\qtnl{110426m}
\frac{m_{x,y}^z}{n}=
\frac{|T|\,c_{\scriptscriptstyle R,\scriptscriptstyle S}^{\scriptscriptstyle T}}
{|R|\,|S|},
\eqtn
where $n=|L|$, $R=\pi(x)$, $S=\pi(y)$, and $T=\pi(z)$. To this end,
choose $g\in\pi(x)$ and count in two ways the set $\Omega$ of all pairs
$(\ell,h)\in L\times S$ such that $g\ell h\in z$.

For each fixed $h\in S$, the number of $\ell\in L$ with $g\ell h\in z$
is $m_{x,y}^z$, since $LgL=x$ and $LhL=y$. Therefore,
\[
|\Omega|=|S|m_{x,y}^z.
\]
On the other hand, set $r=\pi(Lg\ell)\in R$. Since $g\ell=\ell' r$ for
some $\ell'\in L$, the condition $g\ell h\in z$ is equivalent to
$rh\in z\cap G=T$. For each fixed $r\in R$, the number of $\ell\in L$
such that $\pi(Lg\ell)=r$ is $|L|/|R|$. The number of pairs
$(r,h)\in R\times S$ with $rh\in T$ is $|T|c_{R,S}^T$. Consequently,
\[
|\Omega|=\frac{|L|}{|R|}\,|T|c_{R,S}^T.
\]
Comparing the two expressions for $|\Omega|$ yields~\eqref{110426m}.

To complete the proof, note that the bijection $x\mapsto x_{\pi(x)}$,
$x\in X$, is an isomorphism from $\fX$ to $\fX(\cA(K,L))$. Indeed, the
equality in~\eqref{110426f} for the normalized multiplicities follows
from~\eqref{110426m}, and preservation of the involution follows from
$\pi(x^\star)=x^{-1}\cap G=(x\cap G)^{-1}=\pi(x)^{-1}$.
\eprf

Let a finite group $A$ act on a finite group $G$ via a homomorphism
$\varphi:A\to\aut(G)$. In the semidirect product $K=G\rtimes_\varphi A$,
the naturally embedded subgroup~$G$ is a complement to
$L=\{1_G\}\rtimes A$. Hence Theorem~\ref{020426a} gives the following.

\crllrl{110426t}
The coset multivalued group of the triple $(G,A,\varphi)$ is isomorphic
to the multivalued group associated with the S-ring $\fL(G,\varphi(A))$
over $G$. In particular, every coset multivalued group is a Schur
multivalued group.
\ecrllr

\section{A family of multivalued groups that are not double-coset groups}
\label{300326b}
\subsection{Affine S-rings}
We use the S-ring construction of~\cite[Theorem~1]{Hanaki2022}. Let $\mF$
be a finite field of order $q$, let $G=(\mF^2,+)$, and let
$L=\{\lambda I_2:\lambda\in\mF^\times\}$ be the group of scalar matrices.
The affine group $K=G\rtimes L$ acts on $G$, with $G$ acting by
translations and $L$ by scalar multiplication. Consider the S-ring
$\cA_0=\fL(G,L)$ and its S-subrings.

The nonidentity basic sets of $\cA_0$ are the one-dimensional
$\mF$-subspaces of $\mF^2$ with zero removed. There are $q+1$ such sets;
each has size $q-1$ and satisfies $-R=R$. For any two such sets $R,S$,
\qtnl{120426a}
\und{R}\cdot\und{S}=\css
\und{G^\#\setminus(R\cup S)},&\text{if $R\ne S$;}\\[1mm]
(q-1)\,\und{0_G}+(q-2)\und{R},&\text{if $R=S$,}
\ecss
\eqtn
where $G^\#=G\setminus\{0_G\}$.

We construct an S-subring of $\cA_0$ corresponding to a partition
$\pi=(k_1,\ldots,k_a)$ of $q+1$, where
\[
a\ge1,\qquad k_1\ge k_2\ge\ldots\ge k_a\ge1,\qquad
\sum_{i=1}^a k_i=q+1.
\]
Choose a partition $\cP$ of the set $\cS_0^\#$ of nonidentity basic sets
into classes $\cP_u$ of size~$k_u$, $1\le u\le a$. Let $R_u$ be the
union of the sets in $\cP_u$; thus
\qtnl{140426b}
|R_u|=(q-1)k_u,\quad 1\le u\le a.
\eqtn
A direct calculation using~\eqref{120426a} shows that, for
$1\le u,v\le a$,
\qtnl{140426a}
\und{R_u}\und{R_v}=\css
\begin{aligned}
&(k_u-1)k_v\,\und{R_u}+k_u(k_v-1)\und{R_v}\\
&\qquad{}+k_uk_v\,\und{S_{u,v}}
\end{aligned}&\text{if $u\ne v$,}\\[2mm]
\begin{aligned}
&(q-1)k_u\,\und{0_G}+[q-2+(k_u-1)(k_u-2)]\,\und{R_u}\\
&\qquad{}+k_u(k_u-1)\,\und{S_{u,u}}
\end{aligned}&\text{if $u=v$,}
\ecss
\eqtn
where $S_{u,v}=G^\#\setminus(R_u\cup R_v)$. Set $R_0=\{0_G\}$. Then the
$\mC$-vector space
\[
\cA(q,\pi)=\Span_{\mC}\{\und{R_u}:0\le u\le a\}
\]
is a commutative S-ring over $G$. Its basic sets are $R_0,R_1,\ldots,R_a$.
Its structure constants are determined
solely by $k_1,\ldots,k_a$. Hence its algebraic isomorphism class depends
only on $\pi$ and not on the choice of $\cP$. We call $\cA(q,\pi)$ the
\emph{affine} S-ring over a group of order $q^2$ corresponding to the
partition~$\pi$.

\subsection{Construction}
Applying Theorem~\ref{151125a} to affine S-rings yields the following
multivalued groups.

\thrml{120426h}
Let $q$ be a prime power and let $\pi=(k_1,\ldots,k_a)$ be a partition
of~$q+1$. Then the multivalued group $\fX(q,\pi)$ with elements
\[
x_u=\frac{1}{|R_u|}\und{R_u},\qquad 0\le u\le a,
\]
associated with the S-ring $\cA(q,\pi)$ is commutative, has the identity
map as its involution, and has reduced valency
\[
n=(q-1)\cdot\operatorname{LCM}\{k_1,\ldots,k_a\}.
\]
For $1\le i,j\le a$, the multiplicities are given by the following
formulas; multiplication by $x_0$ is determined by the identity axiom:
\[
i\ne j:\qquad m_{x_i,x_j}^{x_u}=\frac{n}{q-1}\cdot\css
0&\text{if $u=0$,}\\
k_u-1&\text{if $u\in\{i,j\}$,}\\
k_u&\text{if $u\not\in\{0,i,j\}$,}
\ecss
\]
\[
i=j:\qquad m_{x_i,x_i}^{x_u}=\frac{n}{(q-1)k_i}\cdot\css
1&\text{if $u=0$,}\\
q-2+(k_i-1)(k_i-2)&\text{if $u=i$,}\\
k_u(k_i-1)&\text{if $u\not\in\{0,i\}$.}
\ecss
\]
For fixed $q$, the multivalued groups $\fX(q,\pi)$ and $\fX(q,\pi')$
are isomorphic if and only if $\pi=\pi'$.
\ethrm
\prf
The first part follows from Theorem~\ref{151125a} and
formulas~\eqref{140426b} and~\eqref{140426a}. Now let
$f:\fX(q,\pi)\to\fX(q,\pi')$ be an isomorphism. Then $a=a'$,
$f(x^{}_0)=x'_0$, and there exists a permutation $\sigma$ of
$\{1,\ldots,a\}$ such that $f(x^{}_i)=x'_{\sigma(i)}$. Preservation of
the normalized multiplicities gives
\[
\frac1{(q-1)k_i}
=\frac{m_{x_i,x_i}^{x_0}}{n}
=\frac{m_{x'_{\sigma(i)},x'_{\sigma(i)}}^{x'_0}}{n'}
=\frac1{(q-1)k'_{\sigma(i)}}.
\]
Hence $k^{}_i=k'_{\sigma(i)}$ for all $i$. Since the parts of both
partitions are listed in nonincreasing order, $\pi=\pi'$. The converse
follows from the explicit formulas for the multiplicities.
\eprf

\rmrkl{050626a}
Under the assumptions and notation of Theorem~\ref{120426h}, suppose
that $q\ge3$ and $k_i\ge2$ for every $1\le i\le a$. The formulas for
the multiplicities imply that $m_{x_i,x_i}^{x_j}>0$ for every
$1\le i\le a$ and every $0\le j\le a$. Thus any subgroup containing
at least one $x_i$, $i\ne0$, contains all elements of $\fX(q,\pi)$.
Consequently, $\fX(q,\pi)$ has no subgroups $Y$ with
$\{x_0\}\subsetneq Y\subsetneq X$.
\ermrk

We show that, among the groups $\fX(q,\pi)$ for which all parts of $\pi$
are at least $2$, the proportion of double-coset groups tends to zero
as $q\to\infty$.

\thrml{180426a}
Let $q\ge3$ be a prime power. Let $\fK(q)$ be a set of representatives
of the isomorphism classes of the groups $\fX(q,\pi)$ for which all
parts of $\pi$ are at least~$2$, and let $\fK_0(q)$ be its subset of
double-coset groups. There exist absolute constants $c_1,c_2>0$ such that
\qtnl{130626a}
|\fK(q)|\ge2^{c_1\sqrt q}\qaq
|\fK_0(q)|\le2^{c_2(\log q)^2}.
\eqtn
In particular, the ratio $\frac{|\fK_0(q)|}{|\fK(q)|}$ tends to $0$
as $q\to\infty$.
\ethrm
\prf
Applying the formula for the number of partitions $\pi$ in which every
part is at least $2$ (see~\cite[p.~500]{SandorMC2006}), we obtain the
following asymptotic formula for the number of multivalued groups
in $\fK(q)$:
\[
|\fK(q)|\sim\frac{C_1}{(q+1)^{3/2}}2^{C_2\sqrt{q+1}}
\]
for some absolute constants $C_1,C_2>0$. This proves the first
inequality in~\eqref{130626a}.

We now bound $|\fK_0(q)|$. By Lemma~\ref{160426h} and the formulas in
Theorem~\ref{120426h}, the permutation group $\ov K$ corresponding to
a double-coset multivalued group $\fX\in\fK_0(q)$ has degree
\[
\sum_{x\in X}\frac1{p_{x,x^\star}^{x_0}}
=1+\sum_{i=1}^a(q-1)k_i=q^2,
\]
where $X$ is the underlying set of $\fX$. By Remark~\ref{050626a} and
Lemma~\ref{160426h}, this permutation group is primitive. Since
conjugate subgroups of $\sym(q^2)$ give rise to isomorphic multivalued
groups, a coarse upper bound for $|\fK_0(q)|$ is the number of conjugacy
classes of all primitive permutation groups of degree~$q^2$.

By~\cite[Theorem I]{PyberS1997}, the number of conjugacy classes of
primitive groups of degree~$d$ is at most $d^{C\mu(d)}$, where $\mu(d)$
is the largest exponent in the prime factorization of $d$ and $C$ is
an absolute constant. If $q=p^f$, then $\mu(q^2)=2f\le2\log_2q$.
Consequently,
\[
|\fK_0(q)|\le(q^2)^{C\mu(q^2)}\le2^{4C(\log_2q)^2}.
\]
This proves the second inequality in~\eqref{130626a}.
\eprf

\section*{Acknowledgments}

The research of A.\,V. Vasil’ev was carried out within the framework of the Sobolev Institute of Mathematics project FWNF-2026-0017.


\end{document}